\documentclass[12pt]{article}

\usepackage{a4wide}
\usepackage{amssymb}
\usepackage{amsfonts}
\usepackage{amsmath}

\date{}

\newtheorem{proposition}{Proposition}[section]
\newtheorem{theorem}[proposition]{Theorem}
\newtheorem{lemma}[proposition]{Lemma}

\newtheorem{corollary}[proposition]{Corollary}

\def\GK{{\rm  GK}\,}
\def\Kdim{{\rm K.dim }\,}

\def\der{\partial }

\def\nFM0{{\nu }_{F,M_0}}
\def\nFN0{{\nu }_{F,N_0}}
\def\nGN0{{\nu }_{G,N_0}}

\def\N0{ {\bf N}_0 }

\def\t{\otimes}

\def\v{\varphi}
\def\ra{\rightarrow}

\def\Xpm{X^{\pm }}

\def\s{\sigma}

\def\l1{{\lambda}_1}

\def\a{\alpha}
\def\a0{ {\alpha }_0}
\def\a1{ {\alpha }_1}

\def\l{\lambda}

\def\nFGM0{{\nu }_{F,G,M_0}}

\def\nFN0{{\nu}_{F,N_0}}

\def\sm{{\sigma}^m}

\def\sm1{{\sigma}^{-1}}

\def\smtp1{{\sigma}^{-t+1}}

\def\S1{S^{-1}}

\def\Xpm1{X^{\pm 1}_1}

\def\sPM1{{\sigma }^{\pm 1}}
\def\sMP1{{\sigma }^{\mp 1 }}

\def\d{\delta}

\def\G{\Gamma}

\def\OO{{\cal O}}
\def\CA{{\cal A}}

\def\CD{{\cal D}}

\def\Ytm1{Y^{t-1}}
\def\Yim1{Y^{i-1}}

\def\CL{{\cal L}}

\def\CF{{\cal F}}

\def\Aut{{\rm Aut}}

\def\Der{{\rm Der }}

\def\dim{{\rm dim }}

\def\ker{ {\rm ker } }

\def\D{ \Delta }
\def\Ev{ {\rm Ev} }

\def\SL2Z{ {\rm SL}_2({\bf Z}) }

\def\CL{{\cal L}}

\def\Gp1{ G^{1 , 1 } }
\def\P11{ P^{-1 , 1 } }
\def\Pp1{ P^{1 , 1 } }

\def\nCLsr{{}^\nu\kern-2pt {\cal L}^{\sigma , \rho  }}
\def\nP{{}^\nu \kern-2pt P}
\def\nL{{}^\nu\kern-2pt L}
\def\nLL{{}^\nu\kern-2pt \Lambda}
\def\nPsr{{}^\nu\kern-2pt P^{\sigma , \rho  }}
\def\nLsr{{}^\nu\kern-2pt L^{\sigma , \rho  }}
\def\nuCL{{}^\nu\kern-2pt  {\cal L}}
\def\nCLsr{{}^\nu\kern-2pt {\cal L}^{\sigma , \rho  }}
\def\nCL1m{{}^\nu\kern-2pt {\cal L}^{-1 , 1  }}
\def\x1nu{x^\frac{1}{\nu}}
\def\xm1nu{x^{-\frac{1}{\nu}}}

\def\trdeg{{\rm tr.deg}}

\def\ra{\rightarrow }

\def\CT{{\cal T}}

\def\nAM0{{\nu }_{{\cal A},M_0}}
\def\nAN0{{\nu }_{{\cal A},N_0}}

\def\Kdim{ {\rm Kdim } }

\def\Der{ {\rm Der }}

\def\vS{\widetilde{S}}
\def\vS1{\widetilde{S}^{-1}}

\def\Gg{\mathfrak{g}}

\def\tor{{\rm tor}}
\begin{document}


\author{V.\  Bavula  }

\title{ Finite generation of division subalgebras and
of the group of eigenvalues for commuting derivations  or
automorphisms of division algebras}

\maketitle
\begin{abstract}
Let $D$ be a division algebra such that $D\t D^o$ is a Noetherian
algebra, then any division subalgebra of $D$ is a {\em finitely
generated} division algebra. Let $\D $ be a finite set of
commuting derivations or automorphisms of the division algebra
$D$, then the group $\Ev (\D )$ of common eigenvalues (i.e. {\em
weights})  is a {\em finitely generated abelian} group. Typical
examples of $D$ are the quotient division algebra ${\rm Frac} (\CD
(X))$ of the ring of differential operators $\CD (X)$ on a smooth
irreducible affine variety $X$ over a field $K$ of characteristic
zero,  and the quotient division algebra ${\rm Frac} (U (\Gg ))$
of the universal enveloping algebra $U(\Gg )$ of a finite
dimensional Lie algebra $\Gg $. It is proved that the algebra of
differential operators $\CD (X)$ is isomorphic to its opposite
algebra $\CD (X)^o$.

\end{abstract}


\section{Introduction}

Throughout this paper, $K$ is a field, $\t=\t_K$. Noetherian means
left and right Noetherian. For a $K$-algebra $A$, $A^o$ denotes
 the {\em opposite} algebra to $A$ (recall that $A^o=A$ as abelian
 groups but the multiplication in $A^o$ is given by the rule: $a*b=ba$),
 and $A^e:=A\t A^o$ is called the {\em enveloping } algebra of
 $A$. The expressions ${}_AM$, $M_A$, and ${}_AM_A$ means that
 $M$ is respectively a left, right $A$-module, and  an $A$-bimodule.
 {\em Finitely generated division algebra} means a division algebra
 which is generated (as a division algebra) by a finite set of
 elements (i.e. $x_1, \ldots , x_n$ is a set of generators for a division
 $K$-algebra $D$ if $D$ is the only division $K$-subalgebra of $D$ that contains
 $x_1, \ldots , x_n$).

For division algebras {\em finite dimensional}  over $K$ there is
a well-developed theory where (commutative) subfields play a
fundamental role. By contrast, if a division algebra is infinite
dimensional little is known about its division subalgebras.

{\it Question. Suppose that $D$ is a finitely generated division
 $K$-algebra, is any division $K$-subalgebra of $D$ finitely
 generated?}

Certainly this is the case when $D$ is a  field. We will see that
the answer is {\em affirmative} for many popular division
algebras. For a similar question about  {\em subfields} ($=$
commutative division $K$-subalgebras), Resco, Small and Wadsworth
give an affirmative answer in \cite{RSW79}: {\em Let $D$ be a
division algebra over a field $K$ such that $ D\t D^o$ is
Noetherian, then every (commutative) subfield of $D$ containing
$K$ is finitely generated.} One of the crucial steps in their
proof is the following result of Vamos \cite{Va78}: {\em Let $L$
be a field extension of $K$. Then $L\t L$ is Noetherian iff $L$ is
a finitely generated  over $K$.} M. Smith \cite{MSmith73} showed
that there is a division algebra $D$ with centre $K$, containing
two maximal subfields whose transcendence degrees are any two
prescribed cardinal numbers.

Let  $\D =\{ \d_1, \ldots , \d_t\}$ be a set of {\em commuting}
$K$-derivations of a division $K$-algebra $D$. The set $\Ev (\D
):=\{ \l =(\l_1, \ldots , \l_t)\in K^t\, | \, \d_i (u)=\l_i u,
i=1, \ldots , t$ for some $0\neq u\in D\}$ of common eigenvalues
is an {\em additive} subgroup of $K^t$, and the $\D$-{\em
eigen-algebra} $D(\D ):=\bigoplus_{\l \in \Ev (\D )}D_\l $ is a
$\Ev (\D )$-graded algebra where $D_\l :=\{ u\in D\, | \,
\d_i(u)=\l_i u, i=1, \ldots , t\}$, $D_\l D_\mu \subseteq D_{\l
+\mu} $ for all $\l , \mu \in \Ev (\D )$, and $0\neq u\in D_\l $
implies $u^{-1}\in D_{-\l }$.

Let  $\D =\{ \d_1, \ldots , \d_t\}$ be a set of {\em commuting}
$K$-automorphisms of a division $K$-algebra $D$, and let
$K^*:=K\backslash \{ 0\}$ be the multiplicative group of the field
$K$. The set $\Ev (\D ):=\{ \l =(\l_1, \ldots , \l_t)\in K^{*t}\,
| \, \d_i (u)=\l_i u, i=1, \ldots , t$ for some $0\neq u\in D\}$
of common eigenvalues is an {\em multiplicative} subgroup of
$K^{*t}$, and the $\D$-{\em eigen-algebra} $D(\D ):=\bigoplus_{\l
\in \Ev (\D )}D_\l $ is a $\Ev (\D )$-graded algebra where $D_\l
:=\{ u\in D\, | \, \d_i(u)=\l_i u, i=1, \ldots , t\}$, $D_\l D_\mu
\subseteq D_{\l \mu} $ for all $\l , \mu \in \Ev (\D )$, and
$0\neq u\in D_\l $ implies $u^{-1}\in D_{\l^{-1} }$.

The first statement of the next result is an extension of the
mentioned above result of Resco-Small-Wadsworth to division
subalgebras (with a short {\em different} proof given in Section
\ref{PRTHM}).

\begin{theorem}\label{DDaccd}
Let $D$ be a division $K$-algebra such that $D\t D$ is a
Noetherian $D$-bimodule, and let $\D =\{ \d_1, \ldots , \d_t\}$ be
either a set of commuting $K$-derivations or commuting
$K$-automorphisms of the  division $K$-algebra $D$.  Then
\begin{enumerate}
\item $D$ satisfies the ascending chain condition on division
$K$-subalgebras, or equivalently, every division $K$-subalgebra of
$D$ is a finitely generated division $K$-algebra.
 \item The group of eigenvalues $\Ev (\D )$ is a finitely generated
 abelian group, and so $\Ev (\D )=\CT \oplus \mathbb{Z}^r$ where $r$ is
 the rank of the group $\Ev (\D )$ and $\CT $ is a finite abelian
 group.
 \item The eigen-algebra $D(\D )$ is a Noetherian domain which
 isomorphic to an iterated skew
 Laurent extension.  In more
 detail, $D_\CT :=\oplus_{\l \in \CT } D_\l $ is a division
 algebra of right and left dimension $|\CT |$ over the division
 algebra $D_0$, $D(\D )$ is isomorphic to the iterated skew
 Laurent extension $D_\CT [x_1, x_1^{-1}; \s_1]\cdots [x_r,
 x_r^{-1}; \s_r]$ with coefficients from the division algebra $D_\CT$.
 \item For each subgroup $F$ of  $\,\Ev (\D )$, $\CF (F):=\oplus_{\l
 \in F}D_\l $ is a Noetherian domain the quotient division algebra
 ${\rm Frac}(\CF (F))$ of which is $\D $-invariant and $\Ev ({\rm
 Frac}(\CF (F)))=F$, any $\D$-eigenvector $v\in {\rm Frac}(\CF
 (F))_\l$ has the form $u^{-1}w$ for some $0\neq u\in D_\mu$,
 $w\in D_{\l +\mu }$, and $\l , \mu \in F$.
\end{enumerate}
\end{theorem}
{\it Remark 1.} $D\t D$ is a {\em Noetherian} $D$-bimodule iff the
 algebra $D\t D^o$ is {\em Noetherian} iff the  algebra
$D\t D^o$ is {\em left} Noetherian iff the  algebra $D\t D^o$ is
{\em right} Noetherian as  it follows from

\begin{equation}\label{D4D}
{}_DD\t D_D\simeq {}_DD\t ({}_{D^o}D^o)\simeq {}_{D\t D^o}D\t D^o,
\; {}_DD\t D_D\simeq D^o_{D^o}\t D_D\simeq D_D\t D^o_{D^o}\simeq
(D\t D^o)_{D\t D^o}.
\end{equation}

{\it Remark 2.} `Finite generation' is built in in the structure
of the eigen-algebra $D(\D )$ in the sense that it is a finitely
generated algebra  over a finitely generated division algebra.

In Section \ref{APLDX}, it is proved that many division algebras
that appear naturally in applications   satisfy the conditions of
Theorem \ref{DDaccd} (Proposition \ref{constg}, Lemma
\ref{itfgd}), {\em eg} ${\rm Frac}(\CD (X))$ (Corollary
\ref{divDXfg}) and ${\rm Frac} (U(\Gg ))$ (Corollary \ref{DUGN}).


\section{Proof of Theorem \ref{DDaccd}}\label{PRTHM}

 Recall that any torsion free
finitely generated abelian group
 is  a free abelian group of finite rank, and vice versa. Any
 finitely generated abelian group $G$ is isomorphic to $\CT \oplus
 \mathbb{Z}^r$ where
 $r:=\dim_\mathbb{Q}(\mathbb{Q}\t_\mathbb{Z}G)$ is the {\em rank}
 of the group $G$, $\CT $ is the {\em torsion} subgroup of $G$,
 that is the subgroup of $G$ that contains all the elements of
 finite order, it is a finite group.

{\bf Proof of Theorem \ref{DDaccd}}.  $1$. Suppose that inside $D$
one can pick a strictly ascending chain of division
$K$-subalgebras $\G_1\subset \G_2\subset \cdots \subset
\G_n\subset \cdots $, we seek a contradiction; this gives a
strictly ascending chain of $D$-sub-bimodules, $K_1\subset
K_2\subset \cdots \subset K_n\subset \cdots $, where $K_n=\ker
(\phi_n)$ where $\phi_n: D\t D\ra D\t_{\G_{n}}D$, $x\t y\ra
x\t_{\G_{n}} y$, (use the fact that $D$ is a free left and right
$\G_n$-module and tensor product commutes with direct sum), a
contradiction. Hence $D$ satisfies the {\em acc} on division
$K$-subalgebras.

$2$, $3$,  and $4$. The proof of two cases are very similar, so we
will treat them simultaneously by making some adjustments to our
notation. So, let $\D =\{ \d_1, \ldots , \d_t\}$ be either a set
of commuting $K$-derivations or a set of commuting
$K$-automorphisms of the division algebra $D$. In the first case,
$\Ev (\D )$ is an {\em additive} subgroup of $K^t$, in the second
case, $\Ev (\D )$ is a {\em multiplicative} subgroup of $K^{*t}$.
In the second case, we still will write the group operation {\em
additively}, i.e. $\l +\mu$ means $\l \mu$, $-\l $ means
$\l^{-1}$, $0$ means $1$. Let $D_0$ be the set of $\D $-{\em
constants}:
 $D_0:=\cap_{i=1}^t\ker_D(\d_i)$, in the case of derivations; and
 $D_0:=\{ d\in D\, | \, \d_i(d)=d, i=1, \ldots , t\}$, in the case
 of automorphisms. In both cases, $D_0$ is a division subalgebra
 of $D$.

 Given a division algebra $\G $, a group $G$, a group homomorphism
 $\v :G\ra \Aut_K(\G )$, and a `2-cocycle' $G\times G\ra \G^*:=\G
 \backslash \{ 0\}$, $(g,h)\mapsto (g,h)$. A {\em generalized
 crossed product} is an algebra $\G *G=\oplus_{g\in G}\G g$  which
 is a free left $\G $-module with multiplication given by the rule
 $$ ag\cdot bh=a \v (g)(b)(g,h)gh, \;\; a,b\in \G, \; g,h\in G.$$
 It follows from $ag=g \v (g)^{-1}(a)$ that $\G g=g\G\simeq
 \G_\G$, and so $\G *G=\oplus_{g\in G}g\G$ is a free right $\G
 $-module. A `2-cocycle' means that the multiplication of the
 generalized crossed product is
 associative. When $G=\mathbb{Z}$ and $(i,j)=1$ for all $i,j\in
 \mathbb{Z}$, we have, so-called, a {\em skew Laurent extension}
 with coefficients from $\G $ denoted $\G [x, x^{-1}; \s]$ where
 $x$ is the group generator $1$ for $\mathbb{Z}$ and $\s =\v
 (1)\in \Aut_K(\G )$. So, the skew Laurent extension generated by
 $\G $ and $x$, $x^{-1}$ subject to the defining relations $x^{\pm
 1}a=\s^{\pm 1}(a)x^{\pm 1}$ for all $a\in \G$. An {\em iterated
 skew Laurent extension } $A_n:= \G [x_1, x_1;\s_1]\cdots [x_n,
 x_n^{-1}; \s_n]$ is defined inductively as $A_n=A_{n-1}[x_n,
 x_n^{-1}; \s_n]$. Since the division algebra $\G $ is a
 Noetherian algebra then the iterated skew Laurent extension $A_n$
 is a Noetherian algebra (1.17, \cite{GW}).

 For each $\l \in \Ev (\D )$, fix $0\neq u_\l \in D_\l $. Then it
 is easy to see that
the $\D$-{\em eigen-algebra}  is a free (left and right)
$D_0$-module: 
\begin{equation}\label{D=dcpz}
\CD =D(\D ):= \bigoplus_{\l \in \Ev (\D )}D_0u_\l =\bigoplus_{\l
\in \Ev (\D )}u_\l D_0, \;\; u_\l u_\mu =(\l , \mu )u_{\l +\mu },
\;  (\l , \mu ):= u_\l u_\mu u_{\l +\mu}^{-1}\in D_0,
\end{equation}
for $a,b\in D_0$: $au_\l bu_\mu=a(u_\l bu^{-1}_\l )u_\l u_\mu =
a(u_\l au^{-1}_\l )(u_\l u_\mu u^{-1}_{\l +\mu })u_{\l +\mu }$ (In
general, this is not a generalized crossed product but if $\Ev (\D
)\simeq \mathbb{Z}^r$ one can choose generators in such a way that
it is).
 Given any finitely generated
 subgroup $F=\CT \oplus (\oplus_{i=1}^s\mathbb{Z}v_i)$ of $\Ev (\D )$
 where $\CT $ is the torsion part of $F$.
The algebra $D_\CT :=\oplus_{\l \in \CT}D_0u_\l =\oplus_{\l \in
\CT} u_\l D_0$ has left and right dimension $|\CT |<\infty$ over
the division algebra $D_0$ where $|\CT |$ is the order of the
group $\CT $.  The map $$l: D_\CT \ra {\rm
End}_{D_\CT}(D_{\CT}),\;  a\mapsto (l_a:x\mapsto ax),$$ is an
algebra isomorphism where ${\rm End}_{D_\CT}(D_{\CT})$ is the
endomorphism algebra of the {\em right} $D_\CT $-module $D_\CT $.
 For each nonzero element $a\in D_\CT $, $l_a$ is a monomorphism,
 hence it is an isomorphism since (the right dimension over $D_0$)
 ${\rm r.dim}_{D_0}(D_\CT )={\rm r.dim}_{D_0}(aD_\CT )=| \CT |<\infty$. Therefore,
 ${\rm End}_{D_\CT}(D_{\CT})$ is  a division algebra, hence so is
 its isomorphic copy $D_\CT $.  Let
 $F'=\oplus_{i=1}^s\mathbb{Z}v_i$, and so  $F=\CT \oplus F'$. The
 subalgebra $\CF' =\CF (F'):= \oplus_{\l \in F'}D_0u_\l$ of $D (\D )$ is
 isomorphic to the iterated skew Laurent extension $\CL :=
 D_0[x_1, x_1^{-1}; \s_1]\cdots [x_s, x_s^{-1}; \s_s]$ where $\s_i
 (d)=u_{v_i}du_{v_i}^{-1}$ $(d\in D_0)$ and $\s_i (x_j)=\l_{ij}x_j$,
 $j<i$, where $\l_{ij}:=u_{v_i}u_{v_j}u_{v_i}^{-1}u_{v_j}^{-1}\in
 D_0$ (via the $K$-algebra epimorphism $\CL \ra \CF' $, $d\mapsto d$,
 $x_i\mapsto u_{v_i}$, where $d\in D_0$). This follows easily from
 a definition of an iterated skew Laurent extension and the facts
 that $D_0$ is a division algebra, $u_{v_1}^{n_1}\cdots
 u_{v_s}^{n_s}\in D_{n_1v_1+\cdots +n_sv_s}$, and
 $F'=\oplus_{i=1}^s\mathbb{Z}v_i$. Then, by a similar reasoning
 (since $D_\CT$ is a division algebra), the algebra $\CF$
 is isomorphic to the iterated skew Laurent extension
 $D_\CT [x_1, x_1^{-1}; \s_1]\cdots [x_s, x_s^{-1}; \s_s]$ where $\s_i
 (d)=u_{v_i}du_{v_i}^{-1}$ $(d\in D_\CT)$ and $\s_i (x_j)=\l_{ij}x_j$,
 $j<i$,  $\l_{ij}$ are as above.

Since $D_\CT$ is a Noetherian algebra so is the algebra $\CF$.
 So,  $\CF
 $ is a Noetherian domain, let ${\rm Frac}(\CF )$ be its quotient
 division algebra, so any element of ${\rm Frac}(\CF )$ is a
 fraction $a^{-1}b$ for some $0\neq a, b\in \CF $. Note that the
 elements $a $ and $b$ are finite sums $\sum a_\l $ and $\sum b_\l
 $ of eigenvectors $a_\l , b_\l \in D_\l $, $\l \in F$. If $0\neq c=a^{-1}b\in
 D_\mu$ for some $\mu \in \Ev (\D )$, then $ac=b\neq 0$ implies that
 $a_\l c=b_\nu$ for some $a_\l\neq 0$ and $b_\nu\neq 0$ such that $\l +\mu =\nu$, and so
 $c=a_\l^{-1}b_\nu$ and $\mu =\nu -\l $. This proves that any
 $\D $-eigenvector of ${\rm Frac}(\CF )$
  is a fraction of the eigenvectors of $\CF $ and that
\begin{equation}\label{EvDFrF}
\Ev (\D , {\rm Frac}(\CF (F) ))=F.
\end{equation}
It follows immediately from this fact and statement 1 that $\Ev
(\D )$ is finitely generated: otherwise one can find in $\Ev (\D
)$ a strictly ascending chain of subgroups: $F_1\subset F_2\subset
\cdots $, which gives, by (\ref{EvDFrF}), the strictly ascending
chain of division subalgebras: ${\rm Frac}(\CF (F_1) )\subset {\rm
Frac}(\CF (F_2) )\cdots $, a contradiction. This finishes the
proof of statement 2 and 4. Then statement 3 has, in fact, been
proved above.  This finishes the proof of
 Theorem \ref{DDaccd}.
 $\Box $

For an abelian monoid $E$, the set $\tor (E)$ of all the elements
$e\in E$ such that $ne=0$ is a group, so-called, the {\em torsion
subgroup} of $E$.

\begin{corollary}\label{1DDaccd}
Let a $K$-algebra $A$ be a Noetherian  domain with $D:={\rm
Frac}(A)$ such that $D\t D$ is a Noetherian $D$-bimodule, $\D =\{
\d_1, \ldots , \d_t\}$ be either a set of commuting
$K$-derivations or commuting $K$-automorphisms of the algebra $A$.
Then the  abelian monoid of eigenvalues $\Ev (\D , A)$ for $\D $
in $A$ is a submonoid of a finitely generated abelian group, and
so the rank of $\Ev (\D , A)$ is finite, and  the torsion subgroup
$\tor (\Ev (\D, A))$ is a finite group.
\end{corollary}

{\it Proof}. Note that each derivation  (resp. an automorphism)
$\d_i$ of $A$ can be uniquely extended to a derivation (resp. an
automorphism) of the division algebra $D$ by the rule $\d_i
(s^{-1}a)=s^{-1}a-s^{-1}\d_i (s)s^{-1}a$ (resp.
$\d_i(s^{-1}a)=\d_i(s)^{-1}\d_i(a)$). So, the zero derivation
$[\d_i, \d_j]=0$ of $A$ has zero extension to $D$, and by
uniqueness it must be zero on $D$. Similarly, the identity
automorphism $[\d_i, \d_j]=\d_i\d_j\d_i^{-1}\d_j^{-1}$ of $A$ has
the obvious  extension to $D$, and by uniqueness it must be the
identity map on $D$. So, $\D$ is a set of commuting derivations
(resp. automorphisms)  of $D$. Clearly, $\Ev (\D , A)\subseteq \Ev
(\D, D)$, and the result follows from Theorem \ref{DDaccd}.(2).
$\Box $

\begin{corollary}\label{2DDaccd}
Let a $K$-algebra $A$ be a commutative affine  domain with
$D:={\rm Frac}(A)$,  $\D =\{ \d_1, \ldots , \d_t\}$ be either a
set of commuting $K$-derivations or commuting $K$-automorphisms of
the algebra $A$. Then the  abelian monoid of eigenvalues $\Ev (\D
, A)$ for $\D $ in $A$ is a submonoid of a finitely generated
abelian group, and so the rank of $\Ev (\D , A)$ is finite, and
the torsion subgroup $\tor (\Ev (\D, A))$ is a finite group.
\end{corollary}

In general, the eigen-algebra  $D(\D )$ is not a finitely
generated algebra even in the case of a commutative affine domain
$A$ since, in general, the $\D$-constants $D_0$ is not a finitely
generated algebra (Hilbert 14'th problem, etc).


\section{Applications}\label{APLDX}

 Let $\G $ be a
$K$-algebra, $\s$ be a $K$-automorphism of $\G $, and $\d \in
\Der_K(\G )$ be a $\s $-{\em derivation} of $\G$: $\d (ab)=\d (a)
b +\s (a) \d (b)$ for $a,b\in \G$. The {\em Ore} extension $A=\G
[x; \s , \d ]$ is a $K$-algebra generated freely by $\G $ and an
element $x$ satisfying the defining relations: $xa=\s (a)x+\d (a)$
for $a\in \G$. Let $\G^o$ be the opposite algebra with
multiplication given by the rule $a*b=ba$. Then $\s \in
\Aut_K(\G^o)$ as $\s (a*b)=\s (ba)=\s (b)\s (a)=\s (a)*\s (b)$,
and so $\s^{-1}\in \Aut_K(\G^o)$, and finally  $\d \s^{-1}\in
\Der_K(\G^o)$ is a $\s^{-1}$-derivation of the algebra $\G^o$:
\begin{eqnarray*}
\d\s^{-1}(a*b)&=&\d\s^{-1}(ba)=\d (\s^{-1}(b)\s^{-1}(a))=
\d\s^{-1}(b)\s^{-1}(a)+b\d\s^{-1}(a)\\
&=&\d\s^{-1}(a)*b+\s^{-1}(a)*\d\s^{-1}(b).
\end{eqnarray*}
\begin{equation}\label{Oreop=Ore}
A^o=\G^o[x;\s^{-1}, -\d\s^{-1}].
\end{equation}
{\it Proof}. The $K$-algebra $A$ is generated by $\G $ and $x$
that satisfy the defining relations:
$x\s^{-1}(a)=ax+\d\s^{-1}(a)$, $a\in \G $, since $\s $ an {\em
automorphism} of $\G $. Hence the $K$-algebra $A^o$ is generated
by the $\G^o$ and $x$ that satisfy the defining relations:
$x*a=\s^{-1}(a)*x-\d\s^{-1}(a)$, $a\in \G^o$, and we are done.
$\Box$

The {\it iterated Ore} $A=\G [x_1; \s_1 , \d_1 ]\cdots [x_n; \s_n
, \d_n ]$ is defined inductively as $$(\G [x_1; \s_1 , \d_1
]\cdots [x_{n-1}; \s_{n-1} , \d_{n-1} ]) [x_n; \s_n , \d_n ].$$ By
(\ref{Oreop=Ore}) and induction on $n$, 
\begin{equation}\label{itOreop}
(\G [x_1; \s_1 , \d_1 ]\cdots [x_n; \s_n , \d_n ])^o\simeq \G^o
[x_1; \s_1^{-1} , -\d_1\s_1^{-1} ]\cdots [x_n; \s_n^{-1} ,
-\d_n\s_n^{-1} ].
\end{equation}
The tensor product of two iterated Ore extensions $A=\G [x_1; \s_1
, \d_1 ]\cdots [x_n; \s_n , \d_n ]$  and $B=\D [y_1; \tau_1 ,
\der_1 ]\cdots [y_m; \tau_m , \der_m ]$ is again an iterated Ore
extension
$$ A\t B=\G \t \D [x_1; \s_1
, \d_1 ]\cdots [x_n; \s_n , \d_n ][y_1; \tau_1 , \der_1 ]\cdots
[y_m; \tau_m , \der_m ]$$ where $\s_i$, $\d_i$ and $\tau_j$,
$\der_j$ act {\em trivially} on the elements where they have not
been defined. Recall that if $\G $ is a domain (resp. a Noetherian
algebra) then so is the iterated Ore extension $A$. If
$\G^e=\G\t\G^o$ is a Noetherian algebra then so is the algebra $\G
$.

\begin{lemma}\label{itfgd}
Suppose that $\G^e=\G\t\G^o$ is a Noetherian domain (then so is
$\G $ and an iterated Ore extension $A=\G [x_1; \s_1 , \d_1
]\cdots [x_n; \s_n , \d_n ]$). Let $D={\rm Frac} (A)$. Then
$D^e=D\t D^o$ is a Noetherian domain, and so the results of
Theorem \ref{DDaccd} hold. \end{lemma}

{\it Proof}. $A^e=A\t A^o$ is an iterated Ore extension with
coefficients from the Noetherian domain  $\G^e$, hence $A^e$ is a
Noetherian domain (1.12, \cite{GW}), and so is its localization
$D^e$. $\Box $

\begin{corollary}\label{divDXfg}
Let $X$ be a smooth irreducible affine variety over a
 field $K$ of characteristic zero, $\CD (X)$
be the ring of differential operators on $X$, and $D(X)={\rm
Frac}(\CD (X))$ be its quotient division algebra. Then $D(X)\t
D(X)^o$ is a Noetherian domain, and so the results of Theorem
\ref{DDaccd} hold.
\end{corollary}

{\it Proof}. The coordinate algebra $\OO=\OO (X)$ of the variety
$X$ is a finitely generated domain with the field of fractions,
say $\G ={\rm Frac}(\OO )$. Let $S=\OO \backslash \{ 0\}$. Then,
by 15.2.6, \cite{MR},  $\S1 \CD (X)\simeq \G [t_1;
\frac{\der}{\der x_1}]\cdots [t_n; \frac{\der}{\der x_n}]$ is an
iterated Ore extension (with trivial automorphisms: $\s_i={\rm
id}_\G $) where $n=\dim (X)$ (the dimension of $X$),  $\G $
contains a rational function field $Q_n:=K(x_1, \ldots, x_n)$
where the
 $\frac{\der}{\der x_i}$ are partial derivations (extended uniquely
from $Q_n $ to $\G$). Note that $\G \t \G^o=\G \t \G$ is a
Noetherian domain as the localization of the domain $\OO(X\times
X)\simeq \OO (X)\t \OO (X)$, the variety $X\times X$ is a smooth
irreducible affine variety. By Lemma \ref{itfgd}, $D(X)^e$ is a
Noetherian domain and every division $K$-subalgebra of $D(X)$ is a
finitely generated division $K$-algebra. This proves the first two
statements.  Then statement 3 follows from Theorem
\ref{DDaccd}.(2). $\Box $

\begin{lemma}\label{AAoFFo}
Let $A$ be a $K$-algebra and $A\ra A^o$, $a\mapsto a^o$, be the
canonical anti-isomorphism $((\l a +\mu b)^o=\l a^o+\mu b^o$ and
$(ab)^o=b^o*a^o$ for all $\l , \mu \in K$ and $a,b\in A$). Then
\begin{enumerate}
\item $(s^{-1})^o=(s^o)^{-1}$ for each unit $s\in A$.
 \item If $S$ is a left (resp. right) Ore subset of $A$ then $S^o$
 is a right (resp. left) Ore subset of $A^o$ and
  $(\S1 A)^o\simeq A^o(S^o)^{-1}$, $s^{-1}a\mapsto a^o* (s^o)^{-1}$
  (resp. $( A\S1 )^o\simeq (S^o)^{-1}A^o$, $as^{-1}\mapsto (s^o)^{-1}*a^o$)
  is the $K$-algebra isomorphism.
  \item  If $A$ is a Noetherian domain then ${\rm
  Frac}(A^o)\simeq {\rm Frac}(A)^o$, $s^{-1}a\mapsto a^o*
  (s^o)^{-1}$, is the isomorphism of division $K$-algebras.
\item  If $A$ is a Noetherian domain such that $A\simeq A^o$ then
 ${\rm  Frac}(A^o)\simeq {\rm Frac}(A)^o$.
\end{enumerate}
\end{lemma}

{\it Proof}. $1$. $ss^{-1}=s^{-1}s=1$ implies
$s^o*(s^{-1})^o=(s^{-1})^o*s^o=1$, and so $(s^{-1})^o=(s^o)^{-1}$.

$2$. Straightforward.

$3$. It is a particular case of statement 2.

$4$. By the universal property of localization, $A\simeq A^o$
implies ${\rm Frac}(A)\simeq {\rm Frac}(A^o)$, and by statement 3,
${\rm Frac}(A^o)\simeq {\rm Frac}(A)^o$. $\Box $

\begin{corollary}\label{DUGN}
Let $\Gg$ be a finite dimensional Lie algebra over a field $K$,
$U=U(\Gg )$ be its universal enveloping algebra, $D(\Gg )={\rm
Frac}(U)$ be its quotient division algebra. Then  $D(\Gg )^e\simeq
D(\Gg )\t D( \Gg )$ is a Noetherian domain, and so the results of
Theorem \ref{DDaccd} hold.
\end{corollary}

{\it Proof}. $U\simeq U^o$, $g\mapsto -g$, $g\in \Gg $. Hence,
${\rm Frac} (U)\simeq {\rm Frac} (U^o)\simeq{\rm Frac} (U)^o$
(Lemma \ref{AAoFFo}) and  $U^e=U\t U^o\simeq U\t U\simeq U(\Gg
\oplus \Gg )$, and so $D(\Gg )^e\simeq D(\Gg ) \t D(\Gg )$ is a
Noetherian domain as a localization of $U(\Gg \oplus \Gg )$, the
rest follows from Theorem \ref{DDaccd}. $\Box $

We have a great similarity in the proofs of the last two
statements, one can repeat this pattern for other `constructions'
of algebras and their division algebras. To formalize the proofs
in many similar situations let us introduce a concept of a {\em
good construction} of  algebras. We say that we have a {\em
construction} of $K$-algebras, say $\CA $, if, for a given
$K$-algebra $\G $, one attaches a set (class) of $K$-algebras $\CA
(\G )$. Examples in mind are  Ore extensions $\CA (\G )=\{ \G
[x;\s, \d ]\}$, iterated Ore extensions, iterated skew Laurent
polynomial algebras, etc. We say that the construction $\CA $ is
{\em good} if the following three properties hold:

 $(G1)$ if $\G $
is a Noetherian domain then so is each algebra from the set $\CA
(\G )$,

$(G2)$ $\CA (\G )^o\subseteq \CA (\G^o )$, and

$(G3)$ $\CA (\G )\t \CA (\G ')\subseteq \CA (\G \t \G')$,

where $\CA (\G )^o:=\{ A^o\, | \, A\in \CA (\G )\}$ and similarly
 $\CA (\G )\t \CA (\G '):=\{ A\t A'\, | \, A\in \CA (\G ), A'\in  \CA ( \G')
 \}$.

For the definitions and properties of the algebras from
 the examples below  the reader is refereed to \cite{GW} and  \cite{MR}.

{\it Examples of good constructions}. $(1)$ Iterated Ore
extensions.

$(2)$ Iterated skew Laurent extensions: $\G [x_1, x_1^{-1}; \s_1]
\cdots [x_n, x_n^{-1}; \s_n]$ where $\s_i$ are $K$-automorphisms
($(G1)$ - use the leading term and  1.17, \cite{GW}; $(G2)$ -
Exercise 1P, p. 17, \cite{GW}; $(G3)$ - obvious).

\begin{proposition}\label{constg}
Suppose that the enveloping algebra  $\G^e=\G \t \G^o$ of an
algebra  $\G $ is a Noetherian domain, $\CA $ is a good
construction, then each algebra  $A\in \CA (\G )$ is a Noetherian
domain and so is $D^e=D\t D^o$ where $D={\rm Frac}(A)$. Hence,
Theorem \ref{DDaccd} is true for $D$.
\end{proposition}

{\it Proof}. $\G^e$ is a  Noetherian domain, then so is $\G $, and
then each algebra $A\in \CA (\G )$ is a Noetherian domain since
the construction $\CA $ is good. $A^e=A\t A^o\in \CA (\G) \t \CA
(\G )^o\subseteq \CA (\G) \t \CA (\G^o )\subseteq \CA (\G \t
\G^o)$, and so $A^e$ is a Noetherian domain, as $\CA$ is good.
Hence, so is its localization $D^e$. The rest is obvious. $\Box $

So, the algebras that satisfy conditions of Theorem \ref{DDaccd}
are fairly common.

\begin{lemma}\label{DNGN}
Let $D$ be a division algebra over a field $K$ such that $D\t D$
is a Noetherian $D$-bimodule. Let $\G $ be a division
$K$-subalgebra of $D$. Then
\begin{enumerate}
\item $\G \t \G$ is a Noetherian $\G $-bimodule, and
 \item $\Kdim ({}_\G(\G \t \G )_\G )\leq \Kdim ({}_D(D\t D)_D)$.
\end{enumerate}
\end{lemma}
{\it Remark}. $\Kdim ({}_\G M_\G)$ means the Krull dimension of a
$\G $-bimodule $M$.

{\it Proof}. Suppose that $\G \t \G$ is not a Noetherian $\G
$-bimodule, we seek a contradiction. Then one can find a strictly
ascending chain of $\G$-sub-bimodules: $I_1\subset I_2\subset
\cdots $. Note that $D\t D\simeq D \t_\G (\G \t \G )\t_\G D$, an
isomorphism of $D$-bimodules. Since $D_\G$ is free, $D\t_\G
I_1\subset D\t_\G I_2\subset \cdots $ is a strictly ascending
chain of $(D, \G )$-bimodules. Similarly, since ${}_\G D$ is free,
$D\t_\G I_1\t_\G D\subset D\t_\G I_2\t_\G D\subset \cdots $ is a
strictly ascending chain of $D$-bimodules, a contradiction. We
have proved that any strictly ascending chain of $\G $-bimodules
$\{ I_i\}$ gives (by tensoring as above) the strictly ascending
chain of $D$-bimodules $\{ D\t_\G I_i\t_\G D\}$, hence
 $\Kdim ({}_\G(\G \t \G )_\G )\leq \Kdim ({}_D(D\t D)_D)$.  $\Box $

Given a $K$-algebra $A$, a $K$-linear map $\s :A\ra A$ is called
an {\em anti-isomorphism} iff $\s (ab)=\s (b) \s (a)$ for all
$a,b\in A$. Clearly, $\s $ is an anti-isomorphism iff $\s : A\ra
A^o$, $a\mapsto \s (a)^o$, is a $K$-algebra isomorphism.

\begin{lemma}\label{AFFDual}
Let $A_i$, $i\in I$, be subalgebras of a $K$-algebra $A$,
$B:=\cap_{i\in I}A_i$, $\s : A\ra A$ be an anti-isomorphism such
that $\s (A_i)=A_i$ for all $i$. Then $\s $ induces the
anti-isomorphism of the algebra $B$.
\end{lemma}

{\it Proof}. Clearly, $\s^{-1}$ is an anti-isomorphism of the
algebra $A$ such that $\s^{-1}(A_i)=A_i$ for all $i\in I$. Then
$\s (B)\subseteq (B)$ and $\s^{-1}(B)\subseteq B$  for all $i\in
I$, hence $\s (B)=B$, and we are done. $\Box $

\begin{theorem}\label{KDXDX}
Let $X$ be a smooth irreducible affine variety over a
 field $K$ of characteristic zero, $\CD (X)$
be the ring of differential operators on $X$, and $D(X)={\rm
Frac}(\CD (X))$ be its quotient division algebra. Then
\begin{enumerate}
\item  $\CD (X)\simeq \CD(X)^o$, $d\mapsto d$, $\der \mapsto -\der
$, where $d\in \OO (X)$ and $\der \in \Der_K(\OO (X))$. \item
$D(X)\simeq D(X)^o$. \item if $X$ is a smooth affine variety (not
necessarily irreducible) over $K$ then still $\CD (X)\simeq \CD
(X)^o$.
\end{enumerate}
\end{theorem}
{\it Remark}. $1$. The algebra $\CD (X)$ is generated by the
coordinate algebra $\OO (X)$ and the $\OO (X)$-module $\Der_K(\OO
(X))$ of $K$-derivations of the algebra $\OO (X)$ (5.6,
\cite{MR}).

$2$. So, the ring of differential operators on a smooth
irreducible algebraic
 variety is {\em symmetric} object indeed. If $A$ is {\em not} smooth
   then, in general, the algebra $\CD
(A)$ need not be  a finitely generated algebra nor a left or right
Noetherian algebra, \cite{BGGDiffcone72}, the algebra $\CD (A)$
can be finitely generated and right Noetherian yet not left
Noetherian, \cite{SmStafDifopcurve}.

 {\it Proof}. $1$. We keep the notation of the proof of Corollary
 \ref{divDXfg}. In particular, the algebra $\CD (X)$ is a
 subalgebra of its localization $A:= \S1 \CD (X)\simeq \G [t_1;
\frac{\der}{\der x_1}]\cdots [t_n; \frac{\der}{\der x_n}]$.

  By 2.13 and 2.6, \cite{MR}, there is a finite set of
 elements of the coordinate algebra $\OO (X)$, say $c_1, \ldots ,
 c_s$, such that the natural inclusion $\CD (X)\ra \prod_{i=1}^s
 \CD (X)_{c_i}$ is a
 {\em faithfully flat} extension where $A_i:=\CD (X)_{c_i}$ is the
 localization of $\CD (X)$ at  the powers of the element $c_i$
 such that $\CD (X)_{c_i}=\OO (X)_{c_i}[t_1;
\frac{\der}{\der x_1}]\cdots [t_n; \frac{\der}{\der x_n}]$. Note
that $A_i\subseteq A$. Hence $\CD (X)=\cap_{i=1}^sA_i$ (let
$B:=\cap_{i=1}^sA_i$ then $\CD (X)\subseteq B\subseteq A_i$ for
each $i$, then the localization of the chain of inclusions above
at $c_i$ gives $A_i:=\CD (X)_{c_i}\subseteq B_{c_i}\subseteq A_i$,
and so, by the faithful flatness, we must have $\CD (X)=B$).
 The map $\s :A\ra A$ given by $d\ra d$ $(d\in \G )$, $t_i\mapsto -t_i$, gives
 an anti-isomorphism of $A$ such that $\s (A_i)=A_i$ for all $i$.
 By Lemma \ref{AFFDual}, $\s $ gives the anti-isomorphism of the
 algebra $\CD (X)$.

$2$. By Lemma \ref{AAoFFo}, $D(X)^o={\rm Frac}(\CD (X))^o\simeq
{\rm Frac}(\CD (X)^o)\simeq {\rm Frac}(\CD (X))=D(X)$.

$3$. Then $X\simeq \prod_{i=1}^sX_i$ is a direct product of smooth
irreducible affine varieties $X_i$ over $K$. Then
$$ \CD (X)=\CD (\prod_{i=1}^sX_i)\simeq \prod_{i=1}^s\CD (X_i)\simeq \prod_{i=1}^s\CD
(X_i)^o \simeq (\prod_{i=1}^s\CD (X_i))^o\simeq \CD (X)^o. \;\;\;
\Box
$$

It is well-known fact that if $C$ is a commutative $K$-subalgebra
of the ring of differential operators $\CD (X)$ (see Theorem
\ref{KDXDX}) then the Gelfand-Kirillov dimension $\GK (C)\leq
n:=\dim (X)$
 (i.e. the transcendence degree ${\rm tr.deg}_K ({\rm Frac}
 (C))\leq n$). It follows from Corollary 3.12, \cite{BavCr}, that
 $\Ev (\D )\simeq \mathbb{Z}^r$ and $r\leq n$, where $\D =\{ \d_1,
 \ldots , \d_t\}$ is a set of commuting {\em locally finite}
 $K$-derivations of the algebra $\CD (X)$. A $K$-derivation $\d$
 of an algebra $A$ is called {\em locally finite} if, for each $a\in A$,
 $\dim_K(\sum_{i\geq 0}K\d^i(a))<\infty$. In a view of Corollary
 3.12, \cite{BavCr} and Theorem
\ref{KDXDX}, the author
 propose the following conjecture.

 {\it Conjecture. If $\D =\{ \d_1,
 \ldots , \d_t\}$ is a set of commuting
 $K$-derivations or $K$-automorphisms  of the division algebra $D(X)={\rm Frac}(\CD
 (X))$ then the rank of the abelian group $\Ev (\D )\leq \dim
 (X)$.}

{\it Question 1. Given a division $K$-algebra $D$ over a field $K$
of characteristic zero such the algebra $D\t D^o$ is Noetherian
and
 a set $\D =\{ \d_1,
 \ldots , \d_t\}$ of commuting
 $K$-derivations or $K$-automorphisms  of the division algebra $D$. Is it true that
  the rank of the abelian group $\Ev (\D )\leq
 \Kdim (D\t D^o)$?}

{\it Remark}.  If, in the question above, $D$ is a (finitely
generated) {\em field} then the result is obviously true (as it
follows from (\ref{D=dcpz}) and  $\GK (D(\D ))\leq \GK (D)$  that
$$\GK (D(\D ))={\trdeg}_K(D_0)+{\rm rank}(\Ev (\D ))\leq \trdeg_K(D)=\Kdim(D\t
D)$$ where $D_0$ is the subfield of $\D $-constants in $D$).

%
%

{\it Question 2. For a singular irreducible affine variety find a
necessary and  sufficient condition that $\CD (X)\simeq \CD
(X)^o$.}

{\it Conjecture 2. Let $X$ be a smooth irreducible affine curve
over an algebraically closed field $K$  of characteristic zero,
and $\d \in \Der_K(\CD (X))$.  Then the eigen-algebra $D(\d )$ is
a finitely generated Noetherian algebra}.

Department of Pure Mathematics

University of Sheffield

Hicks Building

Sheffield S3 7RH

UK

email: v.bavula@sheffield.ac.uk

\end{document}